\documentstyle[12pt, a4wide]{article}
\begin{document}
\newcommand{\qed}{\hfill $\Box$ \\[0.2cm]}
\newcommand{\lug}{\mbox{\boldmath{$\langle \! \langle \;$}}}
\newcommand{\dug}{\mbox{\boldmath{$\; \rangle \! \rangle$}}}
\newcommand{\mk}{\mbox{\boldmath{$k$}}}
\newcommand{\pma}{\mbox{\boldmath{$p$}}}
\newcommand{\mga}{\mbox{\boldmath{$\Gamma$}}}
\newcommand{\mep}{\mbox{\boldmath{${\varepsilon}$}}}
\newcommand{\tv}{\mbox{\rm T}}
\newcommand{\La}{\mbox{$\Lambda$}}
\newcommand{\Lx}{\mbox{$\Lambda_{\times}$}}
\newcommand{\pa}{\mbox{\underline{$a$}}}
\newcommand{\pb}{\mbox{\underline{$b$}}}
\newcommand{\px}{\mbox{\underline{$x$}}}
\newcommand{\mj}{\mbox{\bf 1}}
\newcommand{\ccc}{\mbox{\bf CartCl}}
\newcommand{\frd}{\mbox{\bf Finord}}
\newcommand{\cc}{\mbox{\bf Cart}}
\newcommand{\sym}{\mbox{\bf SyMon}}
\newcommand{\symcl}{\mbox{\bf SyMonCl}}
\newcommand{\skup}{\mbox{\bf Set}}
\newcommand{\dkz}{{\bf Proof:}{\hspace{0.3cm}} \nopagebreak}
\newcommand{\str}{\rightarrow}

\baselineskip=1.05\baselineskip

\title
{The Maximality
of the Typed Lambda Calculus \\ and of Cartesian
Closed Categories}

\author{Kosta Do\v sen and Zoran Petri\' c}

\date{ }
\maketitle
\begin{abstract}
\noindent From the analogue of B\"{o}hm's Theorem proved for the
typed lambda calculus, without product types and with them, it is
inferred that every cartesian closed category that satisfies an equality
between arrows not satisfied in free cartesian closed categories must be a
preorder. A new proof is given here of these results, which were
obtained previously by Richard Statman and Alex K. Simpson.
\end{abstract}
\vspace{0.5cm}
{\em Mathematics Subject Classification (1991)}: 03B40, 18D15, 18A15, 03G30
\vspace{0.5cm}
\section{Introduction}
\noindent In [7] we have shown that every cartesian category
that satisfies an equality between arrows not satisfied in cartesian
categories freely generated by sets of objects must be a preorder; i.e.,
arrows with the same source and target must be equal. In
this paper we give a new proof of the result of [16] (Theorem 1)
that cartesian closed categories are maximal
in the same sense. This means that all equalities between
arrows not assumed for the axiomatization of cartesian, or cartesian
closed categories, are equivalent with each other. Each of them entails
all the other equalities.

It should be stressed that the equalities in question are in the
language of free cartesian categories, or free cartesian closed
categories, and their satisfaction is taken to be
universal with respect to
objects; i.e., atomic symbols for objects are assumed to be variables,
and the equalities are said to be satisfied when they hold for every
assignment of objects to these variables.

The maximality of cartesian categories we have proved previously cannot
be inferred from the maximality of cartesian closed categories we are
working on in this paper, because not every cartesian category need be closed.
So the latter result cannot be simply dubbed a ``generalization'' of
the former. These results are independent, since the inference of the
latter result from the former is impossible too.

The proof of the maximality of cartesian closed categories we are going
to present here is more demanding than the proof of [7]. This
new proof is based on the analogue of B\"{o}hm's Theorem for the typed lambda
calculus, a version of which was established in [19] (Theorem 2).
We prove this analogue first for the typed lambda calculus
with only functional types, in a way different from Statman's,
and from that we pass to the analogue of
B\"{o}hm's Theorem for the typed lambda calculus with product types added.
The proof of the latter analogue reduces essentially to the proof of
the former. These analogues of B\"{o}hm's Theorem cannot be deduced from
the proof of B\"{o}hm's Theorem for the untyped lambda calculus.

To pass from the analogue of B\"{o}hm's Theorem for the typed lambda
calculus with product types to the maximality result for cartesian
closed categories we rely on the categorial equivalence between typed
lambda calculuses and cartesian closed categories, whose discovery is
due to Lambek (see [11], I.11, and references in the
Historical Perspective and Historical Comments on Part I of that
book). This fundamental equivalence, coupled with the understanding
of cartesian closed categories as theories of deduction in
intuitionistic logic, expresses what is usually called the Curry-Howard
correspondence, but which, with more fairness, could be called {\em the
Curry-Lambek-Howard correspondence}. (The Curry-Howard correspondence is
often called an isomorphism, but the term {\em isomorphism} is
more problematic than the looser term {\em correspondence}. If typed lambda
terms are just used as codes for natural-deduction proofs, then there is
presumably an isomorphism between the codes and the things coded, but no
independent algebraic description is given of the things coded. If such an
independent description is given
with the language of cartesian
closed categories, then we fall upon Lambek's {\em equivalence} of
categories, and not upon an {\em isomorphism} of categories.)

The maximality of cartesian and cartesian closed categories is analogous
to the property of the classical propositional calculus
called {\em Post-completeness}. That this calculus is Post-complete
means that if
we add to it a new axiom schema in the language of this calculus, then
we can prove every formula. B\"{o}hm's Theorem in the lambda calculus, or
rather its immediate corollaries, are sometimes termed
``Post-completeness''.

The equational theory of Boolean algebras is also maximal,
i.e. Post-complete. If we add to this theory a new equality
in the language of the theory, then we can deduce
$1=0$ and every other equality. The maximality of cartesian and
cartesian closed categories is analogous to this maximality of
Boolean algebras. Only in equalities we must take care of types,
while in Boolean algebras there is only one type. Another difference is
that in Boolean algebras the equalities producing the extension may
involve variable terms, while in our extensions of the equational
theories of cartesian and cartesian closed categories we envisage
only equalities with constant arrow terms; variables occur only at
the level of types, i.e. the level of objects.

The import of the maximality of cartesian closed categories
for logic is that, in the implication-conjunction fragment of
intuitionistic logic, the choice of equalities between deductions
induced by $\beta\eta$ normalization in natural deduction
is optimal. These equalities, which correspond to the equalities of
cartesian closed categories, are wanted, and no equality is missing,
because any further equality would lead to collapse: all deductions
with the same premises and conclusion would be equal. The import of
the maximality of cartesian categories for conjunctive logic is the same.

Although the results of this paper were already established in
[19] and [16], our proof is different, and
we hope it might shed some new light on the matter. For his proof,
Simpson relies essentially, among other things,
upon a syntactic result of [18] (Theorem 3;
for a proof of this theorem see [14]), which reduces types in
equalities to a particular type, whereas we
rely on a different result from the same paper [18] (Theorem 2),
proved previously in [17] (Theorem 2), which is a
finite-model property for the typed lambda calculus.
Our approach provides an alternative proof of the
type-reducing result of [18] (Theorem 3).

An analogue of B\"{o}hm's Theorem in the typed lambda calculus
without product types is proved in [19] (Theorem 2), without
mentioning B\"{o}hm's Theorem. Statman has even a semantic notion
of consistent extension, rather than a syntactic notion, such as
we have, following B\"{o}hm. (The two notions happen to be
equivalent, however.) Our analogues of B\"{o}hm's Theorem in the
typed lambda calculus, with and without product types, are closer
to standard formulations of this theorem, and our proof is
different from Statman's, which, as Simpson's proof, relies on the
type-reducing result of [18] (Theorem 3). There are, however, some
analogies in the general spirit of these proofs.

The possibility of proceeding as we do is indicated briefly in
[16] (last paragraph of section 5). Simpson says:
``It is an interesting fact that an alternative direct proof
of Theorem 3 is possible using a typed version of the B\"{o}hm-out
technique [1] (Chapter 10). The details are beyond
the scope of this paper.'' (Simpson's Theorem 3 amounts to our
Maximality Corollary in Section 6 below.) We don't know what B\"{o}hm-out
technique Simpson had in mind, but he assured us his approach
is different from ours. Anyway, we couldn't find such a technique
by imitating [1]. Our technique has some intrinsic
difficulties, but presumably not more than the technique of
[19]. Our presentation takes a little bit
more space because we have
tried to help the reader by going into more details. These details,
which were beyond the scope of Simpson's paper, fall exactly
within the scope of ours.

\section{B\"{o}hm's Theorem}

\noindent B\"{o}hm's Theorem in the untyped lambda calculus says that if
$a$ and $b$ are two different lambda terms in $\beta\eta$ normal form,
and $c$ and $d$ are arbitrary lambda terms, then one can
construct terms $h_{1}, \ldots , h_{n}$, $n \geq 0$, and find variables
$x_{1}, \ldots ,x_{m}$, $m \geq 0$, such that
\[ (\lambda_{x_{1} \ldots x_{m}} a) h_{1} \ldots h_{n} = c, \]
\[ (\lambda_{x_{1} \ldots x_{m}} b) h_{1} \ldots h_{n} = d \]
are provable in the $\beta$ lambda calculus (see [1], Chapter 10,
\S4, Theorem 10.4.2; [4], Chapter 11F, \S8, Theorem 5;
[10], Chapitre V, Th\' eor\` eme 2; we know the original paper
of B\"{o}hm [3] only from references). As a corollary of this theorem one
obtains that if $a$ and $b$ are two lambda terms having a normal form
such that $a=b$ is not provable in the $\beta\eta$ lambda calculus
and this calculus is extended with $a=b$, then one can prove every
equality in the extended calculus.

We will demonstrate first the analogue of B\"{o}hm's Theorem in the typed
lambda calculus with only functional types. The standard
proof of B\"{o}hm's Theorem, which may be found in the books cited above,
cannot be transferred to the typed case. At crucial places it introduces
lambda terms that cannot be appropriately typed. For example, for
$\lambda_{xy}xy$ and $\lambda_{xy}x(xy)$ (i.e., the Church numerals for $1$
and $2$) with $x$ of type $p \str p$ and $y$ of type $p$ there is
no appropriate permutator of type $p \str p$ with whose help these two terms
can be  transformed into terms with a head original head normal form
(see [1], Chapter 10, \S3). A more involved example
is given with
the terms $\lambda_{x}x\lambda_{y}(x\lambda_{z}y)$ and
$\lambda_{x}x\lambda_{y}(x\lambda_{z}z)$ with $x$ of type $(p \str p) \str p$
and $y$ and $z$ of type $p$ (we deal with these two typed terms
in the Example of Section 6).

One cannot deduce our analogue of B\"{o}hm's Theorem for the typed lambda
calculus from B\"{o}hm's Theorem for the untyped lambda calculus. The typed
calculus has a more restricted language and does not allow everything
permitted in the untyped case. Conversely, one cannot deduce B\"{o}hm's
Theorem for the untyped lambda calculus from our typed version of this
theorem. Our result covers only cases where $a$ and $b$ are typable by
the same type.

\section{The typed lambda calculus}

\noindent The formulation of the typed lambda calculus with only functional
types we rely on is rather standard (see, for example,
[1], Appendix 1, or [9]). However, we sketch
this formulation briefly, to fix notation and terminology.

{\em Types} are defined inductively by a nonempty set of {\em atomic types}
and the clause ``if $A$ and $B$ are types, then $(A \str B)$ is a type''.
For atomic types we use the schematic letters $p$, $q$, $r$, $\ldots$,
$p_{1}$, $\ldots$, and for all types we use the schematic letters $A$, $B$,
$C$, $\ldots$, $A_{1}$, $\ldots$  We write $A_{B}^{p}$ for the result of
substituting $B$ for $p$ in $A$. ({\em Substitution} means as usual
{\em uniform replacement}.)

{\em Terms} are defined inductively in a standard manner. We have
infinitely many {\em variables} of each type, for which we use the schematic
letters $x$, $y$, $z$, $\ldots$, $x_{1}$, $\ldots$
For arbitrary terms we use the schematic letters
$a$, $b$, $c$, $\ldots$, $a_{1}$, $\ldots$
That a term $a$ is of
type $A$ is expressed by $a : A$. However, for easier reading, we will not
write types inside terms, but will specify the types of variables
separately. For application we use the standard notation, with the standard
omitting of parentheses. For lambda abstraction we will write $\lambda_{x}$
with subscripted $x$, instead of $\lambda x$ (this way we can do without
dots in $\lambda_{x}x$, which is otherwise written $\lambda x . x$).
We abbreviate $\lambda_{x_{1}} \ldots \lambda_{x_{n}} a$ by
$\lambda_{x_{1} \ldots x_{n}} a$, as usual. We write $a_{b}^{x}$
for the result of substituting $b$ for $x$ in $a$, provided $b$ is free
for $x$ in $a$.

If $a$ is a term, let a {\em type-instance} of $a$ be obtained by
substituting some types for the atomic types
in the variables of $a$.

A {\em formula} of the typed lambda calculus $\Lambda$ is of the form
$a=b$ where $a$ and $b$ are terms of the same type.

The calculus $\Lambda$ of $\beta\eta$ equality is axiomatized with the usual
axioms
\[
\begin{array}{lll}
(\beta){\mbox{\hspace{1em}}} &
(\lambda_{x} a) b = a_{b}^{x}, & {\mbox{ provided $b$ is free for $x$ in $a$,}}
\\[0.3cm]
(\eta) &
\lambda_{x} a x = a, & {\mbox{ provided $x$ is not free in $a$,}}
\end{array}
\]
and the axioms and  rules for equality, i.e. $a=a$ and the rule of
replacement of equals. It is not usually noted that the equality of
$\alpha$ conversion can be proved from the remaining axioms as follows:
\begin{eqnarray*}
\lambda_{x} a & = & \lambda_{y} (\lambda_{x} a) y,{\mbox{ by }}(\eta), \\
             & = & \lambda _{y} a_{y}^{x},{\mbox{ by }}(\beta),
\end{eqnarray*}
where $y$ is a variable not occurring in $a$.

\section{Lambda terms for P-functionals}

\noindent Let $P$ be a finite ordinal. In what follows an interesting
$P$ will be greater than or equal to the ordinal $2$. The set of
{\em $P$-types} is defined inductively by specifying that $P$ is a
$P$-type and that if $A$ and $B$ are $P$-types, then $A \str B$, i.e.
the set of all functions with domain $A$ and codomain $B$, is a $P$-type.
Symbols for $P$-types are types with a single atomic type $P$. It is
clear that for $P$ nonempty a $P$-type cannot be named by two different
$P$-type symbols.

An element of a $P$-type is called a {\em $P$-functional}. It is clear
that every $P$-functional is finite (i.e., its graph is a finite set of
ordered pairs) and that in every $P$-type there are only finitely
many $P$-functionals. For $P$-functionals we use the Greek letters
$\varphi$, $\psi$, $\ldots$, $\varphi_{1}$, $\ldots$

Our aim is to define for every $P$-functional a closed term defining it,
in a sense to be made precise. But before that we must introduce a series of
preliminary definitions. In these definitions we take that the calculus
$\Lambda$ is built over types with a single atomic type, which we call
$p$.

Let the type $A_{0}$ be $p$ and let the type $A_{n+1}$ be $A_{n} \str
A_{n}$. For $i \geq 0$, let the type $N_{i}$ be $A_{i+2}$, i.e.
$(A_{i} \str A_{i}) \str (A_{i} \str A_{i})$.

Let $x^{0}(y)$ be $y$ and let $x^{n+1}(y)$ be $x(x^{n}(y))$. The
terms $[n]_{i}$, called {\em Church numerals of type} $N_{i}$, are defined by
\[ [n]_{i} =_{def} \lambda_{xy} x^{n}(y) \]
for $x : A_{i+1}$ and $y : A_{i}$.

For $x$, $y$ and $z$ all of type $N_{i}$, $u : A_{i+1}$, and $v$ and $w$
of type  $A_{i}$, let
\[ C_{i} =_{def} \lambda_{xyzuv} x (\lambda_{w} zuv) (yuv). \]
These are conditional function combinators, because
in the calculus $\Lambda$ one can prove
\[ C_{i} [n]_{i} ab  = \left\{ \begin{array}{ll}
                                a    & \mbox{if $n=0$}\\
                                b    & \mbox{if $n \not= 0$}
                               \end{array}
                       \right. \]

For $x : N_{i+1}$, $y$ and $z$ of type $A_{i+1}$, and $u$ and $v$
of type $A_{i}$, let
\[ R_{i} =_{def} \lambda_{xy} x (\lambda_{zu} y(zu)) (\lambda_{v} v). \]
These combinators reduce the types of numerals; namely,
in $\Lambda$ one can prove
\[ R_{i} [n]_{i+1} = [n]_{i}. \]

For $x$ and $y$ of type $N_{i+1}$, let the exponentiation combinators
be defined by
\[ E_{i} =_{def} \lambda_{xy} x (R_{i} y). \]
In $\Lambda$ one can prove
\[ E_{i} [n]_{i+1} [m]_{i+1} = [m^{n}]_{i}. \]
For $E_{i} a b$ we use the abbreviation $b^{a}$.

For $x$ and $y$ of type $N_{i}$, $z : A_{i+1}$ and $u : A_{i}$, let
the addition and multiplication combinators be defined by
\[
\begin{array}{lll}
S_{i} & =_{def} & \lambda_{xyzu} xz (yzu), \\[0.3cm]
M_{i} & =_{def} & \lambda_{xyzu} x (yz) u.
\end{array}
\]
In $\Lambda$ one can prove
\[
\begin{array}{lll}
S_{i} [n]_{i} [m]_{i} & = & [n+m]_{i}, \\[0.3cm]
M_{i} [n]_{i} [m]_{i} & = & [n \cdot m]_{i}.
\end{array}
\]
For $M_{i} a b$ we use the abbreviation $a{\cdot}b$.

For $x$, $y$ and $z$ of type $N_{i}$, and $u : N_{i+1}$, let the
pairing and projection combinators be defined by
\[
\begin{array}{lll}
\Pi_{i} & =_{def} & \lambda_{xyz} C_{i} zxy, \\[0.3cm]
\pi_{i}^{1}  & =_{def} & \lambda_{u} u [0]_{i}, \\[0.3cm]
\pi_{i}^{2} & =_{def} & \lambda_{u} u [1]_{i}.
\end{array}
\]
In $\Lambda$ one can prove
\[ \pi_{i}^{1} (\Pi_{i} ab) = a, \]
\[ \pi_{i}^{2} (\Pi_{i} ab) = b. \]

For $x : N_{i+1}$ and $y : N_{i+3}$, let
\[
\begin{array}{ll}
T_{i} =_{def} & \lambda_{x} \Pi_{i} (S_{i} [1]_{i} (\pi_{i}^{1} x))
(\pi_{i}^{1} x), \\[0.3cm]
H_{i} =_{def} & \lambda_{y} y T_{i} (\Pi_{i} [0]_{i} [0]_{i}),
\\[0.3cm]
P_{i} =_{def} & \lambda_{y} \pi_{i}^{2}(H_{i} y).
\end{array}
\]
The terms $T_{i}$ and $H_{i}$ are auxiliary,
while the terms $P_{i}$ are predecessor combinators,
because, for $n \geq 1$, one can prove in $\Lambda$
\[
\begin{array}{ll}
P_{i}[n]_{i+3} = & [n-1]_{i}, \\[0.3cm]
P_{i}[0]_{i+3} = & [0]_{i}.
\end{array}
\]
Typed terms corresponding to all the terms $C_{i}$, $R_{i}$, up to
$P_{i}$, may be
found in [2] (cf. [15] and [8]).

For $x$ and $y$ of type $N_{i}$, $z : A_{i+1}$, and $u$ and $v$
of type $A_{i}$, let
\[ Z_{i+1} =_{def} \lambda_{xyzu} x (\lambda_{v} yzu) (zu). \]
These combinators raise the types of numerals for 0 and 1; namely,
in $\Lambda$ one can prove
\[ Z_{i+1} [0]_{i} = [0]_{i+1}, \]
\[ Z_{i+1} [1]_{i} = [1]_{i+1}. \]
The equality $(\eta)$ is essential to prove this.

For $x : N_{i}$, let
\[ D_{i}^{0} =_{def} \lambda_{x} C_{i} x [0]_{i} [1]_{i} \]
and for $k \geq 1$ and $i \geq 3k$ let
\[ D_{i}^{k} =_{def} \lambda_{x} C_{i} x [1]_{i} Z_{i}(Z_{i-1}(Z_{i-2}
(D_{i-3}^{k-1}(P_{i-3}x)))). \]
These combinators check whether a numeral stands for $k$; namely,
for $n \geq 0$, one can prove in $\Lambda$
\[ D_{i}^{k} [n]_{i}  = \left\{ \begin{array}{ll}
                               [0]_{i}    & \mbox{if $n=k$}\\
                               {\mbox{$[1]_{i}$}}    & \mbox{if $n \not= k.$}
                               \end{array}
                       \right. \]

For every $P$-type symbol $A$, let $A^{i}$ be the type obtained from
$A$ by substituting $N_{i}$ for $P$. Now we are ready to define for every
$P$-functional $\varphi \in A$ a closed term $\varphi^{\lambda} : A^{i}$.

Take a $P$-functional $\varphi \in A$, where $A$ is
$B_{1} \str ( \ldots \str (B_{k} \str P) \ldots )$. By induction on the
complexity of the $P$-type symbol $A$ we define a natural
number $\kappa (\varphi)$ and for every $i \geq \kappa (\varphi)$ we define
a term $\varphi^{\lambda} : A^{i}$.

If $A$ is $P$, then $\varphi$ is an ordinal in $P$. Then $\kappa(n) =0$
and $n^{\lambda} : N_{i}$ is $[n]_{i}$ for every $i \geq 0$.

Suppose $k \geq 1$ and $B_{1}$ is $B \str (C \str P)$. It is enough to
consider this case, which gives the gist of the proof. When
$B_{1}$ is $C_{1} \str (C_{2} \str \ldots (C_{l} \str P) \ldots )$
for $l$ different from $2$ we proceed analogously, but with more notational
complications if $l \geq 3$. For $B = \{ \beta_{1}, \ldots ,\beta_{m} \}$ and
$C = \{ \gamma_{1}, \ldots ,\gamma_{r} \}$, by the induction hypotheses,
we have defined  $\kappa(\beta_{1})$, $\ldots$, $\kappa(\beta_{m})$,
$\kappa(\gamma_{1})$, $\ldots$, $\kappa(\gamma_{r})$, for every
$i \geq \kappa(\beta_{1})$ we have defined $\beta_{1}^{\lambda}$, and
analogously for $\beta_{2}$, $\ldots$, $\beta_{m}$, $\gamma_{1}$,
$\ldots$, $\gamma_{r}$. For $B_{1} = \{ \psi_{1}, \ldots ,\psi_{q} \}$,
let $\varphi(\psi_{j}) = \xi_{j} \in B_{2} \str (\ldots \str (B_{k} \str P)
\ldots)$. (Note that $\varphi$ is not necessarily one-one.) By the induction
hypothesis, we have defined $\kappa(\xi_{1}), \ldots , \kappa(\xi_{q})$,
for every $i \geq \kappa(\xi_{1})$ we have defined $\xi _{1}^{\lambda}$,
and analogously for $\xi_{2}, \ldots, \xi_{q}$.

Let now
\[
\begin{array}{llcl}
(\psi_{1}(\beta_{1}))(\gamma_{1})=d_{1} \in P, &
(\psi_{1}(\beta_{2}))(\gamma_{1})=d_{r+1} \in P, &
\ldots &
(\psi_{1}(\beta_{m}))(\gamma_{1})=d_{(m-1)r +1} \in P \\
(\psi_{1}(\beta_{1}))(\gamma_{2})=d_{2} \in P, &
(\psi_{1}(\beta_{2}))(\gamma_{2})=d_{r+2} \in P, &
\ldots &
(\psi_{1}(\beta_{m}))(\gamma_{2})=d_{(m-1)r +2} \in P \\
{\mbox{\hspace{3em}}}\vdots &
{\mbox{\hspace{3em}}}\vdots &
&
{\mbox{\hspace{3em}}}\vdots \\
(\psi_{1}(\beta_{1}))(\gamma_{r})=d_{r} \in P, &
(\psi_{1}(\beta_{2}))(\gamma_{r})=d_{2r} \in P, &
\ldots &
(\psi_{1}(\beta_{m}))(\gamma_{r})=d_{mr} \in P
\end{array}
\]
Let $n_{1} = 2^{d_{1}} \cdot 3^{d_{2}} \cdot \ldots \cdot p_{mr}^{d_{mr}}$,
where $p_{mr}$ is the $mr$-th prime number. Analogously, we obtain the natural
numbers $n_{2}, \ldots , n_{q}$, all different, that correspond to
$\psi_{2}, \ldots , \psi_{q}$.

We can now define $\kappa(\varphi)$ as
\[ max \{ 3 \cdot max \{n_{1}, \ldots , n_{q} \} +1, \kappa(\beta_{1}),
\ldots , \kappa(\beta_{m}), \kappa(\gamma_{1}), \ldots ,
\kappa(\gamma_{r}), \kappa(\xi_{1}), \ldots , \kappa(\xi_{q}) \}.\]
For every $i \geq  \kappa(\varphi)$ and for $x_{1} : B_{1}^{i}$,
let the term $t$ be defined as
\[ [2]_{i}^{\mbox{$x_{1}\beta_{1}^{\lambda}\gamma_{1}^{\lambda}$}} \cdot
[3]_{i}^{\mbox{$x_{1}\beta_{1}^{\lambda}\gamma_{2}^{\lambda}$}} \cdot
\ldots \cdot
[p_{mr}]_{i}^{\mbox{$x_{1}\beta_{m}^{\lambda}\gamma_{r}^{\lambda}$}}
: N_{i-1}. \]

For $x_{2} : B_{2}^{i}, \ldots , x_{k} : B_{k}^{i}$, let
\[
\begin{array}{l}
Q_{1} =_{def}  C_{i} (Z_{i}(D_{i-1}^{n_{1}} t))(\xi_{1}^{\lambda}
x_{2} \ldots x_{k}) Q_{2}, \\[0.3cm]
Q_{2} =_{def}  C_{i} (Z_{i}(D_{i-1}^{n_{2}} t))(\xi_{2}^{\lambda}
x_{2} \ldots x_{k}) Q_{3}, \\[0.3cm]
\vdots \\[0.3cm]
Q_{q-1} =_{def}  C_{i} (Z_{i}(D_{i-1}^{n_{q-1}} t))(\xi_{q-1}^{\lambda}
x_{2} \ldots x_{k}) (\xi_{q}^{\lambda} x_{2} \ldots x_{k}).
\end{array}
\]
We can now, finally, define $\varphi^{\lambda}$ as $\lambda_{x_{1} \ldots
x_{k}} Q_{1}$.

Next we define by induction on the complexity of the $P$-type symbol
$A$, when a $P$-functional $\varphi \in A$ is {\em $i$-defined} by a term
$a :A^{i}$.

We say that a closed term $a : N_{i}$  {\em $i$-defines} an ordinal
$n \in P$ iff in $\Lambda$ we can prove $a = [n]_{i}$.

For a $P$-functional $\varphi \in B \str C$ we say that $a : B^{i} \str C^{i}$
{\em $i$-defines} $\varphi$ iff, for every $\psi \in B$ and every $b : B^{i}$,
if $b$ $i$-defines $\psi$, then $ab : C^{i}$ $i$-defines $\varphi(\psi) \in C$.

We can now prove the following lemma.
\\[0.3cm]
{\bf Lemma 4.1}\hspace{1em}
{\em For every $i \geq \kappa(\varphi)$, the $P$-functional $\varphi \in A$
is $i$-defined by $\varphi^{\lambda} : A^{i}$.}
\\[0.3cm]
\dkz
We proceed by induction on the complexity of the $P$-type symbol $A$.
The case when $A$ is $P$ is trivial.

Let now $A$ be of the form $B_{1} \str ( \ldots \str (B_{k} \str P) \ldots )$
for $k \geq 1$, let  $B_{1} = \{ \psi_{1}, \ldots , \psi_{q} \}$, and
let everything else be as in the inductive step of the definition of
$\varphi^{\lambda}$. Suppose $b_{1} :B_{1}^{i}$ $i$-defines $\psi_{1}$.
We have to check that $\varphi^{\lambda} b_{1}$ $i$-defines
$\varphi(\psi_{1}) = \xi_{1}$.

By the induction hypothesis we have that
$\beta_{1}^{\lambda}$, $\ldots$, $\beta_{m}^{\lambda}$,
$\gamma_{1}^{\lambda}$, $\ldots$, $\gamma_{r}^{\lambda}$,
$\xi_{1}^{\lambda}$, $\ldots$, $\xi_{q}^{\lambda}$
$i$-define $\beta_{1}, \ldots, \beta_{m}, \gamma_{1}, \ldots,
\gamma_{r}, \xi_{1}, \ldots, \xi_{q}$, respectively.
Then we have
\begin{eqnarray*}
\varphi^{\lambda} b_{1} & = & (\lambda_{x_{1} \ldots x_{k}} C_{i} (Z_{i}
(D_{i-1}^{n_{1}} t))(\xi_{1}^{\lambda} x_{2} \ldots x_{k}) Q_{2}) b_{1}
\\
& = &\lambda_{x_{2} \ldots x_{k}} C_{i} (Z_{i}
(D_{i-1}^{n_{1}} t^{x_{1}}_{b_{1}}))
(\xi_{1}^{\lambda} x_{2} \ldots x_{k}) (Q_{2})^{x_{1}}_{b_{1}}.
\end{eqnarray*}
For the closed term $t^{x_{1}}_{b_{1}}$ we have
\[ t^{x_{1}}_{b_{1}} = [2]_{i}^{\mbox{$b_{1}\beta_{1}^{\lambda}
\gamma_{1}^{\lambda}$}} \cdot
[3]_{i}^{\mbox{$b_{1}\beta_{1}^{\lambda}\gamma_{2}^{\lambda}$}} \cdot
\ldots \cdot
[p_{mr}]_{i}^{\mbox{$b_{1}\beta_{m}^{\lambda}\gamma_{r}^{\lambda}$}} . \]
It follows by the induction hypothesis that
$b_{1}\beta_{1}^{\lambda}\gamma_{1}^{\lambda}$ $i$-defines $d_{1}$,
which means that in $\Lambda$ we can prove
$b_{1}\beta_{1}^{\lambda}\gamma_{1}^{\lambda} = [d_{1}]_{i}$. We proceed
analogously with the other exponents. So in $\Lambda$
we can prove $t_{b_{1}}^{x_{1}} = [n_{1}]_{i-1}$. Hence in $\Lambda$
we have  $D_{i-1}^{n_{1}} t^{x_{1}}_{b_{1}} = [0]_{i-1}$,
and we conclude that
\begin{eqnarray*}
\varphi^{\lambda} b_{1} & = & \lambda_{x_{2} \ldots x_{k}} \xi^{\lambda}_{1}
x_{2} \ldots x_{k}
\\
& = & \xi_{1}^{\lambda},{\mbox{ by }}(\eta).
\end{eqnarray*}
So $\varphi^{\lambda}b_{1}$ $ i$-defines $\xi_{1}$.

Suppose now $b_{2} : B_{1}^{i}$ $i$-defines $\psi_{2}$. Then in $\Lambda$
we have
\[
\varphi^{\lambda} b_{2} =  \lambda_{x_{2} \ldots x_{k}} C_{i} (Z_{i}
(D_{i-1}^{n_{1}} t^{x_{1}}_{b_{2}}))
(\xi_{1}^{\lambda} x_{2} \ldots x_{k})
(C_{i} (Z_{i} (D_{i-1}^{n_{2}} t^{x_{1}}_{b_{2}}))
(\xi_{2}^{\lambda} x_{2} \ldots x_{k}) (Q_{3})^{x_{1}}_{b_{2}}).
\]
Since in $\Lambda$ we can prove
$t_{b_{2}}^{x_{1}} = [n_{2}]_{i-1}$, we can also prove
$D_{i-1}^{n_{1}} t^{x_{1}}_{b_{2}} = [1]_{i-1}$,
and we conclude that
\[
\varphi^{\lambda} b_{2} =  \lambda_{x_{2} \ldots x_{k}}
C_{i} (Z_{i} (D_{i-1}^{n_{2}} [n_{2}]_{i-1}))
(\xi_{2}^{\lambda} x_{2} \ldots x_{k}) (Q_{3})^{x_{1}}_{b_{2}}.
\]
Finally, we obtain as above that $\varphi^{\lambda} b_{2}$ $i$-defines
$\xi_{2}$. We proceed analogously for $\psi_{3}, \ldots, \psi_{q}$.
\qed

This lemma does not mean that we can $i$-define all $P$-functionals
simultaneously for some $i$. But we can always find such an $i$ for
finitely many $P$-functionals.

\section{P-models}

\noindent A model based on $P = \{ 0, \ldots, h-1 \}$, with  $h \geq 2$,
for the calculus $\Lambda$ built over types with a single atomic
type $p$ will be defined as in [9].

An {\em assignment} is a function $f$ assigning to a variable $x : A$ of
$\Lambda$ a functional $f(x)$ in the $P$-type $A^{p}_{P}$, where
$A^{p}_{P}$ is obtained from $A$ by substituting $P$ for $p$.
For an assignment $f$ and a variable $y$, the assignment $f^{y}_{\alpha}$
is defined by
\[ f_{\alpha}^{y} (x)  = \left\{ \begin{array}{ll}
                                \alpha    & \mbox{if $x$ is $y$}\\
                                f(x)      & \mbox{if $x$ is not $y$.}
                               \end{array}
                       \right. \]

If $F$ is the set of all $P$-functionals, then
the $P$-{\em model} is a pair $\langle F , V \rangle$ such that $V$
maps the pairs $(a,f)$, with $a$ a term and $f$ an assignment, into $F$.
We write $V_{a,f}$ instead of $V(a,f)$.
The function $V$ must satisfy the conditions
\[
\begin{array}{ll}
& V_{x,f} = f(x),{\mbox{\hspace{11em}}} \\[0.3cm]
& V_{ab,f} = V_{a,f}(V_{b,f}), \\[0.3cm]
{\mbox{for $x : A$ and $\alpha : A^{p}_{P}$}}, &
V_{\lambda_{x} a , f} (\alpha) = V_{a , f^{x}_{\alpha}}.
\end{array}
\]
There is exactly one such function $V$.

Let $a : A$ be a term such that $x_{1} : A_{1}, \ldots, x_{n} : A_{n}$
are all the variables, both free and bound, occurring in $a$. Let $f$ be  an
assignment, and for every $j \in \{ 1, \ldots, n \}$ let $b_{j}$ $i$-define
$f(x_{j})$. Finally, let \pa\ be the type-instance of $a$ obtained by
substituting $N_{i}$ for $p$. The type of \pa\ is $(A^{p}_{P})^{i}$.
Then we can prove the following lemma.
\\[0.3cm]
{\bf Lemma 5.1}\hspace{1em}
{\em  The term $\pa^{\mbox{$\px_{1} \ldots \px_{n}$}}
_{\mbox{$b_{1} \ldots b_{n}$}}$
$i$-defines $V_{a,f}$. }
\\[0.3cm]
The proof proceeds by a straightforward induction on the complexity of the term
$a$.

Of course, when $a$ is closed, $V_{a,f}$ does not depend on $f$,
and has the same value for all assignments $f$. So, for closed terms
$a$, we can write $V_{a}$ instead of $V_{a,f}$, and we shall do so from
now on.

As an immediate corollary of Lemma 5.1 we obtain the following lemma.
\\[0.3cm]
{\bf Lemma 5.2}\hspace{1em}
{\em
If $a$ is closed, then \pa\ $i$-defines $V_{a}$. }

\section{The maximality of the typed lambda calculus}

\noindent We are now ready to prove our analogue of B\"{o}hm's Theorem for the
typed lambda calculus $\Lambda$, which is not necessarily built over types
with a single atomic type.
\\[0.3cm]
{\bf Theorem 6.1}\hspace{1em}
{\em
If $a$ and $b$ are of the same type and $a = b$ is not provable in $\Lambda$,
then for every two terms $c$ and $d$ of the same type one can construct
type-instances $a'$ and $b'$ of $a$ and $b$, respectively, and terms
$h_{1}, \ldots, h_{n}$, $n \geq 0$, and also find variables $x_{1}, \ldots,
x_{m}$, $m \geq 0$, such that
\[ (\lambda_{x_{1} \ldots x_{m}} a') h_{1} \ldots h_{n} = c, \]
\[ (\lambda_{x_{1} \ldots x_{m}} b') h_{1} \ldots h_{n} = d \]
are provable in $\Lambda$.}
\\[0.3cm]
\dkz
Let $a_{1}$ and $b_{1}$ be type-instances of $a$ and $b$, respectively,
obtained by substituting $p$ for all atomic types. It is easy to
see that $a = b$ is provable in $\Lambda$ iff $a_{1} = b_{1}$ is
provable in $\Lambda$.

Let $x_{1}, \ldots, x_{m}$ be all the free variables in $a_{1}$ or $b_{1}$.
Then since $a_{1} = b_{1}$ is not provable in $\Lambda$, the equality
$\lambda_{x_{1} \ldots x_{m}} a_{1} = \lambda_{x_{1} \ldots x_{m}} b_{1}$
is not provable in $\Lambda$. Let $a_{2}$ be
$\lambda_{x_{1} \ldots x_{m}} a_{1}$
and let $b_{2}$ be
$\lambda_{x_{1} \ldots x_{m}} b_{1}$.

It follows from a theorem of [17] (Theorem 2, p. 187)
and [18] (Theorem 2, p. 21) that if $a_{2} = b_{2}$
is not provable in $\Lambda$, then there exists a $P$-model
$\langle F , V \rangle$ such that $V_{a_{2}} \not= V_{b_{2}}$.
Soloviev's and Statman's theorem doesn't mention exactly
$P$-models, which are based on the full type structure
built over an ordinal $P$, but instead it mentions completely
analogous models based on the full type structure built over
a finite set $S$.

We can always name the elements of $S$ by ordinals so that
$S$ becomes an ordinal $P$. Moreover, for every two distinct elements $s_{1}$
and $s_{2}$ of $S$ we can always name the elements of $S$ so that $s_{1}$
is named by $0$ and $s_{2}$  is named by $1$. This means that the elements
of $S$ can always be named by elements of $P$ so that in the $P$-model
$\langle F , V \rangle$ above there are $P$-functionals $\varphi_{1}, \ldots,
\varphi_{k}$, $k \geq 0$, such that
\[
((V_{a_{2}}(\varphi_{1}))(\varphi_{2})) \ldots (\varphi_{k}) = 0,\]
\[
((V_{b_{2}}(\varphi_{1}))(\varphi_{2})) \ldots (\varphi_{k}) = 1.
\]
Take an even $i \geq max \{ \kappa(\varphi_{1}), \ldots, \kappa(\varphi_{k}) \}$.
By Lemma 4.1, the closed terms
$\varphi_{1}^{\lambda}, \ldots, \varphi_{k}^{\lambda}$ $i$-define
$\varphi_{1}, \ldots, \varphi_{k}$, respectively. By Lemma 5.2, the term
$\pa_{2}$ $i$-defines $V_{a_{2}}$ and $\pb_{2}$ $i$-defines $V_{b_{2}}$.
It follows that in $\Lambda$ we can prove
$\pa_{2} \varphi_{1}^{\lambda} \ldots \varphi_{k}^{\lambda} = [0]_{i}$
and
$\pb_{2} \varphi_{1}^{\lambda} \ldots \varphi_{k}^{\lambda} = [1]_{i}$.

For $x : A_{i}$, $y : A_{i-1}$ and $z : A_{i-2}$ we can prove in $\Lambda$
\[ [0]_{i} (\lambda_{xyz} yz) (\lambda_{yz} z) = [0]_{i-2}, \]
\[ [1]_{i} (\lambda_{xyz} yz) (\lambda_{yz} z) = [1]_{i-2}. \]
So there are closed terms $c_{1}, \ldots, c_{i}$ such that in
$\Lambda$ we can prove
\[ \pa_{2} \varphi_{1}^{\lambda} \ldots \varphi_{k}^{\lambda}
c_{1} \ldots c_{i} = [0]_{0}, \]
\[ \pb_{2} \varphi_{1}^{\lambda} \ldots \varphi_{k}^{\lambda}
c_{1} \ldots c_{i} = [1]_{0}. \]
Let the left-hand sides of these two equalities be $a_{3}$ and
$b_{3}$, respectively.

Take now $c$ and $d$ of type $A$ and take the type-instances $a_{4}$
and $b_{4}$  of $a_{3}$ and $b_{3}$, respectively, obtained by
substituting $A$ for $p$. For $u : A$ we can prove in $\Lambda$
\[ a_{4} (\lambda_{u} d) c = c , \]
\[ b_{4} (\lambda_{u} d) c = d . \]
The terms $a_{4}$ and $b_{4}$ are of the form
$(\lambda_{x_{1} \ldots x_{n}} a') h_{1} \ldots h_{k+i}$
and
$(\lambda_{x_{1} \ldots x_{n}} b') h_{1} \ldots h_{k+i}$ .
If $(N_{i})^{p}_{A}$ is obtained by substituting $A$ for $p$
in $N_{i}$, then $a'$ is a type-instance of $a$ obtained by
substituting $(N_{i})^{p}_{A}$ for every atomic type.
\qed

 Since the procedure for constructing $h_{1}, \ldots , h_{n}$ in the
 proof of Theorem 6.1 can be pretty involved, it may be useful to
 illustrate this procedure with an example. For this example we
 take two terms unequal in \La\ that we mentioned in Section 2
 (this is the more involved of the examples given there).
\\[0.3cm]
 {\bf Example:}{\hspace{1em}}
 Let $a$ and $b$ be $\lambda_{x} x \lambda_{y} (x
 \lambda_{z} y)$ and  $\lambda_{x} x \lambda_{y} (x \lambda_{z} z)$,
 respectively, with
 \linebreak
 $x : (p \str p) \str p$, $y : p$ and $z : p$.
 Since all the atomic types of $a$ and $b$ are already $p$, and
 since these two terms are closed, we have that $a_{2}$ is $a$ and
 $b_{2}$ is $b$.

 The $P$-model falsifying $a = b$ has $P = \{ 0,1 \}$ and $P \str P =
 \{ \psi_{1}, \psi_{2}, \psi_{3}, \psi_{4} \}$, where
\[
 \begin{array}{l}
 \psi_{1} (0) = \psi_{1} (1) = 0, \\
 \psi_{2} (0) = \psi_{2} (1) = 1, \\
 \psi_{3} (0) = 0,{\mbox{\hspace{1em}}} \psi_{3} (1) = 1, \\
 \psi_{4} (0) = 1,{\mbox{\hspace{1em}}} \psi_{4} (1) = 0. \\
 \end{array}
 \]
For $\varphi \in (P \str P) \str P$ defined by
\[
\varphi ( \psi_{1} ) = 1,{\mbox{\hspace{1em}}}
\varphi ( \psi_{2} ) = \varphi ( \psi_{3} ) =
\varphi ( \psi_{4} ) = 0
\]
we have $V_{a} ( \varphi ) = 0$ and $V_{b} ( \varphi ) = 1$.

Then
\[
\begin{array}{ll}
n_{1} = 2^{0} \cdot 3^{0} = 1 & {\mbox{corresponds to $\psi_{1}$}},\\
n_{2} = 2^{1} \cdot 3^{1} = 6 & {\mbox{corresponds to $\psi_{2}$}},\\
n_{3} = 2^{0} \cdot 3^{1} = 3 & {\mbox{corresponds to $\psi_{3}$}},\\
n_{4} = 2^{1} \cdot 3^{0} = 2 & {\mbox{corresponds to $\psi_{4}$}},
\end{array}
\]
and $\kappa ( \varphi ) = 19$. For every $i \geq 19$ and for
$x_{1} : N_{i} \str N_{i}$, the term $t$ is defined as
\linebreak
$[2]_{i}^{\mbox{$x_{1} [0]_{i}$}} \cdot
[3]_{i}^{\mbox{$x_{1}[1]_{i}$}} : N_{i-1}$.
The term $\varphi^{\lambda}$ is defined as
\[
\lambda_{x_{1}}
C_{i} (Z_{i}
(D_{i-1}^{1} t)) [1]_{i}
(
C_{i} (Z_{i}
(D_{i-1}^{6} t)) [0]_{i}
(
C_{i} (Z_{i}
(D_{i-1}^{3} t)) [0]_{i} [0]_{i})).
\]

The terms \pa\ and \pb\ are like $a$ and $b$ with $x : (N_{20} \str
N_{20}) \str N_{20}$, $y : N_{20}$ and $z : N_{20}$, and let $i$ in
$\varphi^{\lambda}$ be 20. Then in \La\ we can prove
$\pa \varphi ^{\lambda} = [0]_{20}$ and
$\pb \varphi ^{\lambda} = [1]_{20}$.
The remaining steps in the construction of $a_{3}$ and $b_{3}$
are straightforward, and we shall not pursue this example further.
\\[0.2cm]

By taking that for $x$ and $y$ of the same type the term $c$ is
$\lambda_{xy} x$ and $d$ is $\lambda_{xy} y$, we obtain the following
refinement of Theorem 6.1.
\\[0.3cm]
{\bf Theorem 6.2}\hspace{1em}
{\em
If $a$ and $b$ are of the same type and $a = b$ is not provable in $\Lambda$,
then for every two terms $e$ and $f$ of the same type one can construct
type-instances $a'$ and $b'$ of $a$ and $b$, respectively, and closed terms
$h_{1}, \ldots, h_{l}$, $l \geq 0$, and also find variables $x_{1}, \ldots,
x_{m}$, $m \geq 0$, such that
\[ (\lambda_{x_{1} \ldots x_{m}} a') h_{1} \ldots h_{l} e f = e, \]
\[ (\lambda_{x_{1} \ldots x_{m}} b') h_{1} \ldots h_{l} e f = f \]
are provable in $\Lambda$.}
\\[0.3cm]

It is clear that if $a$ and $b$ are closed, we need not
mention in this theorem the variables $x_{1}, \ldots , x_{m}$
and we can omit the $\lambda$-abstraction $\lambda_{x_{1} \ldots x_{m}}$.

Although our proof of Theorem 6.1 relies on the equality $(\eta)$ at
some key steps (as we noted in connection with the combinator
$Z_{i+1}$), it is possible to derive a strengthening of this
theorem, as well as of Theorem 6.2, where $\Lambda$ is replaced by
$\Lambda_{\beta}$, which is $\Lambda$ minus $(\eta)$ and plus the
equality of $\alpha$ conversion. We learned how to obtain
this strengthening from Alex Simpson.

First note that if a term $a$ is in both contracted and expanded
$\beta \eta$ normal form, and $a = b$ in $\Lambda$, then
$a = b$ in $\Lambda_{\beta}$. For if $a = b$ in $\Lambda$, then,
since $a$ is in contracted $\beta \eta$ normal form, there is
a term $a'$ such that $b$ $\beta$-reduces to $a'$ and $a'$
$\eta$-reduces by contractions to $a$. But then, since $a$
is also in expanded $\beta \eta$ normal form, $a'$ must be the same term
as $a$. So $a = b$ in $\Lambda_{\beta}$.

Then, as we did to derive Theorem 6.2, take in Theorem 6.1 that
$c$ is $\lambda_{xy} x$ and $d$ is $\lambda_{xy} y$ for $x$ and
$y$ of atomic type $p$. The terms $c$ and $d$ are then in
both contracted and expanded $\beta \eta$ normal form, and hence
it is easy to infer Simpson's strengthening mentioned above
by instantiating $p$ with an arbitrary type.

To formulate below a corollary of Theorem 6.1 we must explain what it
means to extend $\Lambda$ with a new axiom. Let $a$ and $b$ be
of type $A$, and let $a'$ and $b'$ be type-instances of $a$ and $b$
respectively. Then assuming $a = b$ as a new axiom in $\Lambda$
means assuming also $a' = b'$. In other words, $a = b$ is assumed as an axiom
schema, atomic types being understood as schematic letters.
The postulate $(\beta)$ and $(\eta)$ are also assumed as axiom
schemata, in the same sense. We could as well add to $\Lambda$
a new rule of substitution for atomic types. The calculus $\Lambda$ is
closed under this substitution rule
(i.e., this rule is admissible, though not derivable from
the other rules). And any extension of $\Lambda$ we envisage should be
closed under this rule. The rule of substitution of types says that
atomic types are variables.

We can now state the following corollary of Theorem 6.1.
\\[0.3cm]
{\bf Maximality Corollary}\hspace{1em}
{\em  If $a = b$ is not provable in $\Lambda$, then in $\Lambda$
extended with $a = b$ we can prove every formula $c = d$.}

\section{The typed lambda calculus with product types}

\noindent We want to demonstrate next the analogue of B\"{o}hm's Theorem in
the typed  $\beta\eta$ lambda calculus with product types, i.e. with
surjective pairing, projections and a constant of terminal type, by
reducing it to our analogue of B\"{o}hm's Theorem for the typed
lambda calculus \La. The idea of this reduction is inspired by
[6] (Chapter 4.1), [17] (pp. 180ff) and [20].

In the typed lambda calculus with product types, types include an atomic
constant type \tv\ and the type-forming operation $\times$ besides
$\str$. Terms now include an atomic constant term $k : \tv$. Moreover,
for every term $a : A \times B$ we have the terms $p^{1}a : A$
and $p^{2}a : B$, and for all terms $a : A$ and $b : B$ we have the term
$\langle a,b \rangle : A \times B$.

The typed lambda calculus \Lx\ of $\beta\eta$ equality is
axiomatized with the postulates for \La\ extended with the axioms
\[
\begin{array}{lll}
(\times \beta) & p^{1} \langle a,b \rangle = a, &
p^{2} \langle a,b \rangle = b,
\\[0.3cm]
(\times \eta) & \langle p^{1}c , p^{2}c \rangle = c, &
\\[0.3cm]
(\tv) & {\mbox{for $x : \tv$}}, & x=k.
\end{array}
\]

\section{Product normal form of types of \Lx}

\noindent Consider the following reductions of types, which consist in
replacing subtypes of the form on the left-hand side by subtypes of the
form on the right-hand side:
\[
\begin{array}{ccc}
{\mbox{\em redexes}} & {\mbox{\hspace{3em}}} & {\mbox{\em contracta}} \\[0.3cm]
A \str (B_{1} \times B_{2}) &     & (A \str B_{1}) \times (A \str B_{2}) \\
(A_{1} \times A_{2}) \str B &     & A_{1} \str (A_{2} \str B)  \\
A \times (B \times C)    &        & (A \times B) \times C \\[0.3cm]
A \str \tv        &               &  \tv  \\
\tv \str B        &               &   B   \\[0.3cm]
A \times \tv      &               &   A    \\
\tv \times A      &               &   A
\end{array}
\]

A type of \Lx\ is in {\em product normal form} iff it does not have
subtypes that are redexes. If $\times^{1}_{i=1} A_{i}$ is $A_{1}$
and $\times^{n+1}_{i=1} A_{i}$ is $(\times^{n}_{i=1}A_{i}) \times A_{n+1}$,
then a type $A$ of \Lx\ in product normal form is either of the form
$\times^{n}_{i=1} A_{i}$ with every $A_{i}$ a type of \La\ or $A$
is simply \tv.

Every type $A$ of \Lx\ can be reduced by the reductions above to a unique
product normal form $A^{\pi}$. This follows from the Church-Rosser property
of these reductions. These reductions are also strongly normalizing.
This can be proved by assigning uniformly to atomic types, including \tv,
natural numbers greater than or equal to $2$, which will be their {\em
complexity measure}. The complexity measure $c(C)$ of a type $C$ is
computed according to
\[ c(A \times B) = (c(A) + 1) c(B), \]
\[ c(A \str B) = c(B)^{c(A)}. \]
Then all reductions decrease the complexity measure. (The
measure of complexity of
[6] Chapter 4.1, Proposition 4.1.2, does not work, and
we couldn't extract a suitable measure from [17].)

A term $a : A \str B$ is an {\em isomorphism} iff there is a term
$a' : B \str A$ such that for $x : A$ and $y : B$ we can prove
$a' (ax) = x$ and $a (a'y) = y$. We have the following lemma.
\\[0.3cm]
{\bf Lemma 8.1}\hspace{1em}
{\em
For every type $A$ one can construct an isomorphism $h : A \str A^{\pi}$.}
\\[0.2cm]

\noindent For the proof see [6] (Chapter 1.9, Theorem 1.9.9).

Let $a^{\nu}$ be the expanded $\beta\eta$ normal form
of a term $a$ of \Lx\ (sometimes also called {\em long normal form};
see [21], Chapter 6.5, and [6],
Chapter 2). This normal form is unique, according to [5]
(see also references in [6]).

Let $\Pi_{i=1}^{1} a_{i}$ be $a_{1}$ and let $\Pi_{i=1}^{n+1} a_{i}$
be $\langle \Pi_{i=1}^{n} a_{i} , a_{n+1} \rangle$. We can easily prove
the following lemmata.
\\[0.3cm]
{\bf Lemma 8.2}\hspace{1em}
{\em
Let $a$ be a closed term of type $A^{\pi}$ for some $A$. Then
$a^{\nu}$  is either of the form $\Pi_{i=1}^{n} a_{i}$ with
$a_{i}$ a term of \La, or $a^{\nu}$ is $k$.}
\\[0.2cm]

\noindent For that we rely on the fact that the type of a subterm of $a^{\nu}$
must be a subtype of $A^{\pi}$, which in logic is called the
{\em subformula property}.
\\[0.3cm]
{\bf Lemma 8.3}\hspace{1em}
{\em
For $a$ and $b$ of type $A$ and $h : A \str A^{\pi}$ an isomorphism,
if $a=b$ is not provable in \Lx, then
\begin{enumerate}
\item $(ha)^{\nu}$ is of the form $\Pi_{i=1}^{n} a_{i}$ for
$a_{i}$ a term of \La,
\item $(hb)^{\nu}$ is of the form $\Pi_{i=1}^{m} b_{i}$ for
$b_{i}$ a term of \La,
\item $n = m$,
\item one can find an $i \in \{ 1, \ldots, n \}$ such that $a_{i} = b_{i}$
is not provable in \La.
\end{enumerate}}
\noindent \dkz
Clauses 1 and 2 follow by Lemma 8.2, since $a$ and $b$ are not of
type \tv.
Otherwise, $a=b=k$ is provable in \Lx. Clause 3
follows from the fact that $ha$ and $hb$ are of the same type.

For clause 4 we have that if $a_{i} = b_{i}$ is provable in \La, then
$a_{i} = b_{i}$ is provable in \Lx. So if for every $i \in \{1, \ldots, n \}$
we had $a_{i} = b_{i}$ provable in \La, then $(ha)^{\nu} = (hb)^{\nu}$
would be provable in \Lx, which would imply that $a = b$ is provable
in \Lx.
\qed

\section{The maximality of \Lx}
\noindent
For $a : \times_{i=1}^{n} A_{i}$ we define $\pi^{i} a : A_{i}$ as $a$ for
$n = 1$, and for $n > 1$,
\[ \pi^{i} a  = \left\{ \begin{array}{ll}
                                p^{2} a    & \mbox{if $i=n$}\\
                                \pi^{i} p^{1} a    & \mbox{if $i <  n$.}
                               \end{array}
                       \right. \]

We can now state the analogue of B\"{o}hm's Theorem for \Lx.
\\[0.3cm]
{\bf Theorem 9.1}\hspace{1em}
{\em
If $a$ and $b$ are of the same type and $a = b$ is not provable in \Lx,
then for every two terms $c$ and $d$ of the same type of \Lx\
one can construct
type-instances $a'$ and $b'$ of $a$ and $b$, respectively, and terms
$h, h_{1}, \ldots, h_{n}$, $n \geq 0$, and also find variables $x_{1}, \ldots,
x_{m}$, $m \geq 0$, and a natural number $i$ such that
{\samepage
\[ \pi^{i}(h \lambda_{x_{1} \ldots x_{m}} a') h_{1} \ldots h_{n} = c, \]
\[ \pi^{i}(h \lambda_{x_{1} \ldots x_{m}} b') h_{1} \ldots h_{n} = d \]}
are provable in \Lx.}
\\[0.3cm]
\noindent The proof is obtained by applying Lemma 8.3 and Theorem 6.1.

It would be possible to prove Theorem 9.1 directly, without passing
through Theorem 6.1, by applying Soloviev's version of
Soloviev-Statman's theorem, which is given for \Lx\ (see [17],
Theorem 2, p. 187). This would yield a different form for the equalities of
Theorem 9.1.

There is an analogue of Theorem 6.2 obtained by refining Theorem 9.1,
and for closed terms $a$ and $b$ of \Lx, the variables
$x_{1}, \ldots, x_{m}$ are not mentioned. We can,
of course, prove an analogue of Theorem 9.1 for
the typed lambda calculus with nonempty product types, i.e. without
terminal type \tv. Finally, we can draw from Theorem 9.1
the Maximality Corollary where \La\ is replaced by \Lx.

In the last part of this work we shall apply the following version of
Theorem 9.1.
\\[0.3cm]
{\bf Theorem 9.2}\hspace{1em}
{\em
If $a$ and $b$ are closed terms of the same type and
$a = b$ is not provable in \Lx,
then one can construct
type-instances $a'$ and $b'$ of $a$ and $b$, respectively, and closed terms
$h, h_{1}, \ldots, h_{l}$, $l \geq 0$,
and also find  a natural number $i$ such that}

{\em \samepage
\[ \pi^{i}(h a') h_{1} \ldots h_{l} (p^{1}x) (p^{2}x) = p^{1}x, \]
\[ \pi^{i}(h b') h_{1} \ldots h_{l} (p^{1}x) (p^{2}x) = p^{2}x \]}
{\em are provable in \Lx.}

\section{Free cartesian closed categories}

\noindent The equational calculus $CCC$ of cartesian closed categories
is introduced as follows. {\em Object terms} of $CCC$ are the types of
\Lx. For {\em arrow terms} of $CCC$ we use the schematic letters
$f$, $g$, $h$, $\ldots$, $f_{1}$, $\ldots$, and we indicate by
$f : A \vdash B$, where $A$ and $B$ are types, that $A$ is the
{\em source} and $B$ the {\em target} of $f$ (we use $\vdash$
instead of the usual $\str$, which we have reserved for a type operation).
We say that $A \vdash B$ is the {\em arrow type} of $f$.

Arrow terms are defined inductively starting from the
{\em atomic arrow terms}

\[ \mj_{A} : A \vdash A ,\]
\[ \pma^{1}_{A,B} :A \times B \vdash A, {\mbox{\hspace{1em}}}
\pma^{2}_{A,B} :A \times B \vdash B, \]
\[\mep_{A,B} : (A \str B) \times A \vdash B, \]
\[\mk_{A} : A \vdash \tv , \]
with the help of the partial operations on arrows
\[ \frac{f : A \vdash B {\mbox{\hspace{3em}}} g : B \vdash C}
{g \circ f : A \vdash C} \]

\[\frac{f : C \vdash A {\mbox{\hspace{3em}}} g : C \vdash B}
{\lug f,g \dug : C \vdash A \times B} \]

\[\frac{f : C \times A \vdash B}
{\mga_{C,A} f : C \vdash A \str B} \]

The {\em formulae} of $CCC$ are equalities $f = g$ where the arrow terms
$f$ and $g$ have the same arrow type.

The axioms of $CCC$ are
\[ f = f, \]
\[ f \circ \mj_{A} = \mj_{B} \circ f = f \; , {\mbox{\hspace{1em}}}
h \circ (g \circ f) = (h \circ g) \circ f, \]
\[ \pma^{1}_{A,B} \circ \lug f,g \dug = f \; , {\mbox{\hspace{1em}}}
\pma^{2}_{A,B} \circ \lug f,g \dug = g, \]
\[ \lug \pma^{1}_{A,B} \circ h \; , \; \pma^{2}_{A,B} \circ h \dug = h, \]
\[ \mep_{A,B} \circ \lug \mga_{C,A} f \circ \pma^{1}_{A \str B , A} \;
, \; \pma^{2}_{A \str B , A} \dug = f, \]
\[ \mga_{C,A} (\mep_{A,B} \circ \lug g \circ \pma^{1}_{C,A} \; , \;
\pma^{2}_{C,A} \dug ) = g, \]
\[ {\mbox{for }} f : A \vdash \tv \; ,  {\mbox{\hspace{1em}}} f = \mk_{A}, \]
and its inference rules are replacement of equals and substitution
of types for atomic types. (Substitution of types says that atomic types
are variables; cf. the remarks on substitution of types at the end of
Section 6.)

We speak of {\em the} calculus $CCC$, but, as a matter of fact, there
are many such calculuses obtained by varying the generating set
of atomic types. When we say that a cartesian closed category $\cal K$
is a model of $CCC$, we mean that the arrow terms $\mj_{A}$ of $CCC$
are interpreted by the unit arrows of $\cal K$, that the arrow terms
$\pma^{1}_{A,B}$ and $\pma^{2}_{A,B}$ of $CCC$ are interpreted by the
projection arrows of $\cal K$, etc.

The calculus $CCC$ does not have arrow-term variables, and hence its models
are not necessarily cartesian closed categories. They need not
even be categories. This calculus catches only the
``canonical-arrow fragment'' of cartesian closed categories.

The term model of $CCC$ is a free cartesian closed category. This
category is isomorphic to the category ${\cal C}(\Lx)$, engendered
by \Lx, which we will define by slightly varying the approach of
[11] (I.11; see also [12]).
The objects of ${\cal C}(\Lx)$
are again the types of \Lx, and the arrows are equivalence classes
$ [ \lambda_{x} a ] = \{ b : A \str B \; \mid \;
b {\mbox{ is a closed term of }}\Lx {\mbox{ and }} b = \lambda_{x} a
{\mbox{ is provable in }} \Lx \}$. The arrow type of $ [ \lambda_{x} a ] $
is $A \vdash B$.

The $CCC$ structure of ${\cal C}(\Lx)$ is defined by
\[ \mj_{A} = [\lambda_{x} x], {\mbox{ for }} x : A, \]
\[\pma^{1}_{A,B} = [ \lambda_{x} p^{1}x ]\; , {\mbox{\hspace{1em}}}
\pma^{2}_{A,B} = [ \lambda_{x} p^{2}x ]\; , {\mbox{ for }}
x : A \times B, \]
\[ \mep_{A,B} =  [ \lambda_{x} p^{1}x(p^{2}x) ]\; , {\mbox{ for }}
x : (A \str B) \times A, \]
\[ \mk_{A} = [ \lambda_{x} k ]\; , {\mbox{ for }} x : A, \]
\[ [a] \circ [b] = [ \lambda_{x} a (bx)], \]
\[ \lug [a] , [b] \dug = [ \lambda_{x} \langle ax , bx \rangle ], \]
\[ \mga_{C,A} [a] = [ \lambda_{xy} a \langle x , y \rangle ]. \]

\section{The maximality of cartesian closed categories}

\noindent We can now prove the following theorem.
\\[0.3cm]
{\bf Maximality of CCC}\hspace{1em}
{\em
If the formula $f = g$ of $CCC$ is not provable in $CCC$, then every
cartesian closed category that is a model of $CCC$ extended
with $f=g$ is a preorder.}
\\[0.3cm]
\dkz
Suppose that for $f,g : A \vdash B$ the equality $f=g$ is not provable in
$CCC$. So one can construct closed terms $\lambda_{x} a$ and
$\lambda_{x} b$ of \Lx\ of type $A \str B$, with $x : A$, such that we have
$[\lambda_{x} a] = f$ and $[\lambda_{x} b] = g$, and
$\lambda_{x} a = \lambda_{x} b$ is not provable in \Lx. By
Theorem 9.2, one can construct type-instances $\lambda_{x'} a'$ and
$\lambda_{x'} b'$ of $\lambda_{x} a$ and $\lambda_{x} b$,
respectively, and closed terms $h, h_{1}, \ldots, h_{l}$, $l \geq 0$,
such that for $y : p \times p$,
\[ \lambda_{y}(\lambda_{z}\pi^{i}(h z) h_{1} \ldots h_{l} (p^{1}y) (p^{2}y))
\lambda_{x'} a' = \lambda_{y} p^{1}y, \]
\[ \lambda_{y}(\lambda_{z}\pi^{i}(h z) h_{1} \ldots h_{l} (p^{1}y) (p^{2}y))
\lambda_{x'} b' = \lambda_{y} p^{2}y \]
are provable in \Lx.

Let $\cal K$ be a cartesian closed category that is a model of
$CCC$ plus $f=g$. Then there is a cartesian closed functor $F$ from
${\cal C}(\Lx)$ to $\cal K$ such that $F(p)=C$, for $C$ an arbitrary
object of $\cal K$. If $\pi^{i}(h z) h_{1} \ldots h_{l} (p^{1}y) (p^{2}y)$
is abbreviated by $c$, then we have in $\cal K$
\[ F([\lambda_{y}(\lambda_{z} c) \lambda_{x'} a']) =
\mep_{F(A') \str F(B'),C} \circ \lug F([\lambda_{yz}c]) \; , \;
\mga_{C \times C, F(A')} (F([\lambda_{x'} a']) \circ \pma^{2}_{C \times C,
F(A')}) \dug \]
and the analogous equality obtained by replacing $a'$ by $b'$.
Since in $\cal K$ we have $F([\lambda_{x'} a']) = F([\lambda_{x'} b'])$,
we obtain in $\cal K$
\[ F([\lambda_{y} p^{1}y]) = F([\lambda_{y} p^{2}y]), \]
i.e., $\pma^{1}_{C,C} = \pma^{2}_{C,C}$. Then for
$h_{1}, h_{2} : E \vdash C$ in $\cal K$ we have
\[ \pma^{1}_{C,C} \circ \lug h_{1} \; , \; h_{2} \dug =
   \pma^{2}_{C,C} \circ \lug h_{1} \; , \; h_{2} \dug, \]
i.e., $h_{1} = h_{2}$.
\qed

To prove the Maximality of CCC one could also use
the Maximality Corollary for \Lx\
and the soundness of \Lx\ with respect to $CCC$ models
(see [13], Theorem 4.2, p. 310).

Note that the construction of $a'$, $b'$, $h$, $h_{1}$, $\ldots$, $h_{n}$
and $i$ in Theorem 9.2, as well as in our other analogues of
B\"{o}hm's Theorem, is in principle effective, though
it may be quite involved.
(This relies on the effectiveness of Soloviev's and Statman's theorem
mentioned in the proof of Theorem 6.1.) The construction of
$\lambda_{x} a$ and $\lambda_{x}b$ out of $f$ and $g$ in the proof
of the Maximality of CCC is also effective, and derivations in \Lx\ made
of equalities between closed terms, on which we can rely in the proof of
the Maximality of CCC, can be transformed effectively into derivations
in $CCC$.
All this entails that we have a constructive method to derive
$h_{1} = h_{2}$ from $f = g$ in the proof of the Maximality of CCC.

The Maximality of CCC makes it possible to generalize Soloviev's
and Statman's theorem
we have used in the proof of Theorem 6.1.
(\v Cubri\' c has in [5] a related theorem about the existence
of a faithful cartesian closed functor from free cartesian closed
categories with free arrows, i.e. from models of $CCC$ extended with
arrow-term variables, into the category of sets.) This generalization
says that every equality not provable in $CCC$ can be falsified in any
cartesian closed category that is not a preorder. Soloviev and Statman
envisage as a falsifying category only the category of finite sets. Our
previous result of [7] yields an analogous statement for
cartesian categories. All these matters about generalizing Soloviev's
and Statman's theorem (as well as \v Cubri\' c's) are treated
clearly, systematically and with much insight in [16].
\\[0.5cm]
\begin{center}
{\bf References}
\end{center}
\vskip 10pt

{\baselineskip=0.8\baselineskip
\noindent\hangindent=\parindent
[1] H.P. Barendregt, {\em The Lambda Calculus: Its Syntax
and Semantics}, North-Holland, Amsterdam, 1981, revised edition 1984.

\vskip 5pt

\noindent\hangindent=\parindent
[2] E. Barendsen, {\em Representation of Logic, Data Types and
Recursive Functions in Typed Lambda Calculi}, Doctoraal Scripte,
Faculteit Wiskunde en Informatica, Katholicke Universiteit Nijmegen,
1989.

\vskip 5pt

\noindent\hangindent=\parindent
[3] C. B\"{o}hm, {\em Alcune propriet\` a delle forme $\beta$-$\eta$-normali
nel $\lambda$-$K$-calcolo}, Pubblicazioni dell'Isti\-tuto per le
Applicazioni del Calcolo, Rome, {\bf 696} (1968), 19 pp.

\vskip 5pt

\noindent\hangindent=\parindent
[4] H.B. Curry, J.R. Hindley and J.P. Seldin,
{\em Combinatory Logic, Volume {\rm II}}, North-Holland, Amsterdam, 1972.

\vskip 5pt

\noindent\hangindent=\parindent
[5] D. \v Cubri\' c, {\em Embedding of a free cartesian closed category
into the category of sets}, J. Pure Appl. Algebra
{\bf 126} (1998), 121-147.

\vskip 5pt

\noindent\hangindent=\parindent
[6] R. Di Cosmo, {\em Isomorphism of Types: From $\lambda$-Calculus
to Information Retrieval and Language Design}, Birkha\" user, Boston,
1995.

\vskip 5pt

\noindent\hangindent=\parindent
[7] K. Do\v sen and Z. Petri\' c, {\em The maximality of cartesian
categories}, preprint, Rapport IRIT 97-42 (1997).

\vskip 5pt

\noindent\hangindent=\parindent
[8] S. Fortune, D. Leivant and M. O'Donnel, {\em The expressiveness
of simple and second-order type structures}, J. ACM {\bf 30}
(1983), 151-185.

\vskip 5pt

\noindent\hangindent=\parindent
[9] H. Friedman, {\em Equality between functionals}, in:
R. Parikh ed., {\em Logic Colloquium '73}, Lecture Notes in Math.
{\bf 453}, Springer, Berlin, 1975, 22-37.

\vskip 5pt

\noindent\hangindent=\parindent
[10] J.-L. Krivine, {\em Lambda-calcul: Types et mod\` eles}, Masson,
Paris, 1990 (English translation, Ellis Horwood, 1993).

\vskip 5pt

\noindent\hangindent=\parindent
[11] J. Lambek and P.J. Scott, {\em Introduction to Higher-Order
Categorical Logic}, Cambridge University Press, Cambridge, 1986.
\vskip 5pt

\noindent\hangindent=\parindent
[12] G.E. Mints, {\em Category theory and proof theory} (in Russian),
in: {\em Aktual'nye voprosy logiki i metodologii nauki}, Naukova Dumka,
Kiev, 1980, 252-278 (English translation, with permuted
title, in: G.E. Mints, {\em Selected Papers in Proof Theory}, Bibliopolis,
Naples, 1992).

\vskip 5pt

\noindent\hangindent=\parindent
[13] J.C. Mitchell and P.J. Scott, {\em Typed lambda models and cartesian
closed categories}, in: J.W. Gray and A. Scedrov eds,
{\em Categories in Computer Science and Logic},
Contemp. Math. {\bf 92}, American Mathematical Society,
Providence, 1989, 301-316.

\vskip 5pt

\noindent\hangindent=\parindent
[14] J.G. Riecke, {\em Statman's 1-Section Theorem},
Inform. and Comput. {\bf 116} (1995), 294-303.

\vskip 5pt

\noindent\hangindent=\parindent
[15] H. Schwichtenberg, {\em Definierbare Funktionen im Lambda-Kalk\"{u}l
mit Typen}, Arch. math. Logik Grundlagenforsch. {\bf 17} (1976), 113-114.
(We know this paper only from references.)

\vskip 5pt

\noindent\hangindent=\parindent
[16] A.K. Simpson, {\em Categorical completeness results for the
simply-typed lambda-calculus}, in:
M. Dezani-Ciancaglini and G. Plotkin eds,
{\em Typed Lambda Calculi and Applications (Edinburgh, 1995)},
Lecture Notes in Comput. Sci.
{\bf 902}, Springer, Berlin, 1995, 414-427.

\vskip 5pt

\noindent\hangindent=\parindent
[17] S.V. Soloviev, {\em The category of finite sets and cartesian
closed categories} (in Russian), Zapiski nauchn. sem. LOMI
{\bf 105} (1981), 174-194 (English translation in
J. Soviet Math. {\bf 22}, 1983, 1387-1400).

\vskip 5pt

\noindent\hangindent=\parindent
[18] R. Statman, {\em Completeness, invariance and $\lambda$-definability},
J. Symbolic Logic {\bf 47} (1982), 17-26.

\vskip 5pt

\noindent\hangindent=\parindent
[19] R. Statman, {\em $\lambda$-definable functionals
and $\beta\eta$-conversion},
Arch. math. Logik Grundlagenforsch. {\bf 23} (1983), 21-26.

\vskip 5pt

\noindent\hangindent=\parindent
[20] A.S. Troelstra, {\em Strong normalization for typed terms
with surjective pairing}, Notre Dame J. Formal Logic
{\bf 27} (1986), 547-550.

\vskip 5pt

\noindent\hangindent=\parindent
[21] A.S. Troelstra and H. Schwichtenberg, {\em Basic Proof Theory},
Cambridge University Press, Cambridge, 1996.

\vskip 7pt}

\newpage

\noindent {\em Acknowledgement.} We would like to thank Alex Simpson
for reading a previous version of this paper,
and for making a very helpful suggestion (noted in Section 6).
We are also grateful to Slobodan Vujo\v sevi\' c for his
careful reading of the text and for his comments.

\vskip 5pt

\noindent University of Toulouse III, IRIT, 31062 Toulouse cedex, France
and Mathematical Institute, Knez Mihailova 35, P.O. Box 367, 11001
Belgrade, Yugoslavia, email: kosta@mi.sanu.ac.yu

\vskip 5pt

\noindent University of Belgrade, Faculty of Mining and Geology,
Dju\v sina 7, 11000 Belgrade, Yugoslavia, email: zpetric@rgf.rgf.bg.ac.yu

\vspace{2ex}

\begin{center}
{\bf Note added in September 2012}
\end{center}
\vskip 10pt

\noindent In Theorem~1 of another paper of Richard Statman ({\em
Simply typed $\lambda$ calculus with surjective pairing}, in: H.\
Barendregt et al. eds., {\em Dirk van Dalen Festschrift},
Quaestiones Infinitae, Publications of the Department of
Philosophy, Utrecht University, vol. V, 1993, pp.\ 143-159) one
can find stated what is essentially the maximality result of this
paper.

\end{document}